\documentclass[12pt]{article}
\usepackage{amsmath}
\usepackage{amscd}
\usepackage{amssymb}
\usepackage{array}

\newtheorem{definition}{Definition}[section]
\newtheorem{theorem}{Theorem}[section]
\newtheorem{example}{Example}[section]
\newtheorem{observation}{Observation}[section]
\newtheorem{lemma}{Lemma}[section]

\newtheorem{corollary}{Corollary}[section]
\begin{document}

\newcommand{\id}{\relax{\rm 1\kern-.28em 1}}

\newcommand{\R}{\mathbb{R}}
\newcommand{\C}{\mathbb{C}}
\newcommand{\Z}{\mathbb{Z}}

\newcommand{\g}{\mathfrak{G}}
\newcommand{\e}{\epsilon}

\newcommand{\bp}{\mathbf{p}}
\newcommand{\bmax}{\mathbf{m}}
\newcommand{\m}{\mathbf{m}}

\newcommand{\cA}{\mathcal{A}}
\newcommand{\cB}{\mathcal{B}}
\newcommand{\cC}{\mathcal{C}}
\newcommand{\cI}{\mathcal{I}}
\newcommand{\cO}{\mathcal{O}}
\newcommand{\cG}{\mathcal{G}}
\newcommand{\cJ}{\mathcal{J}}
\newcommand{\cF}{\mathcal{F}}
\newcommand{\cP}{\mathcal{P}}
\newcommand{\ep}{\mathcal{E}}
\newcommand{\E}{\mathcal{E}}
\newcommand{\cl}{\ell}

\newcommand{\rGL}{\mathrm{GL}}
\newcommand{\rSL}{\mathrm{SL}}
\newcommand{\rOSp}{\mathrm{OSp}}
\newcommand{\rsl}{\mathrm{sl}}
\newcommand{\rM}{\mathrm{M}}
\newcommand{\M}{\mathrm{M}}
\newcommand{\End}{\mathrm{End}}
\newcommand{\diag}{\mathrm{diag}}

\newcommand{\fsl}{\mathfrak{sl}}
\newcommand{\fg}{\mathfrak{g}}
\newcommand{\ff}{\mathfrak{f}}
\newcommand{\fgl}{\mathfrak{gl}}
\newcommand{\fosp}{\mathfrak{osp}}

\newcommand{\str}{\mathrm{str}}
\newcommand{\tr}{\mathrm{tr}}
\newcommand{\defi}{\mathrm{def}}
\newcommand{\Ber}{\mathrm{Ber}}
\newcommand{\spec}{\mathrm{Spec}}
\newcommand{\sschemes}{\mathrm{(sschemes)}}
\newcommand{\sschemeaff}{\mathrm{ {( {sschemes}_{\mathrm{aff}} )} }}

\newcommand{\sets}{\mathrm{(sets)}}
\newcommand{\Top}{\mathrm{Top}}
\newcommand{\sarf}{ \mathrm{ {( {salg}_{rf} )} }}
\newcommand{\arf}{\mathrm{ {( {alg}_{rf} )} }}

\newcommand{\odd}{\mathrm{odd}}

\newcommand{\alg}{\mathrm{(alg)}}
\newcommand{\sa}{\mathrm{(salg)}}
\newcommand{\SA}{\mathrm{(salg)}}
\newcommand{\salg}{\mathrm{(salg)}}
\newcommand{\varaff}{ \mathrm{ {( {var}_{\mathrm{aff}} )} } }
\newcommand{\svaraff}{\mathrm{ {( {svar}_{\mathrm{aff}} )}  }}
\newcommand{\ad}{\mathrm{ad}}
\newcommand{\Ad}{\mathrm{Ad}}
\newcommand{\pol}{\mathrm{Pol}}
\newcommand{\Lie}{\mathrm{Lie}}

\newcommand{\al}{\alpha}
\newcommand{\be}{\beta}

\newcommand{\lra}{\longrightarrow}
\newcommand{\ra}{\rightarrow}

\vskip 2cm

    \centerline{\LARGE \bf  On Algebraic Supergroups, Coadjoint Orbits}

  \smallskip

  \centerline{\LARGE \bf and their Deformations }\bigskip

\vskip 2.5cm

\centerline{ R. Fioresi \footnote{Investigation supported by the
University of Bologna, funds for selected research topics.}}

\smallskip

\centerline{\it Dipartimento di Matematica, Universit\`a di
Bologna }
 \centerline{\it Piazza di Porta S. Donato, 5.}
 \centerline{\it 40126 Bologna. Italy.}
\centerline{{\footnotesize e-mail: fioresi@dm.UniBo.it}}

\medskip

\centerline{and} \medskip

\centerline{ M. A. Lled\'o }

\smallskip

 \centerline{\it INFN, Sezione di Torino,
Italy  and} \centerline{\it Dipartimento di Fisica, Politecnico di
Torino,} \centerline{\it Corso Duca degli Abruzzi 24, I-10129
Torino, Italy.} \centerline{{\footnotesize e-mail:
lledo@athena.polito.it}}

 \vskip 3cm

\begin{abstract}

In this paper we study algebraic supergroups and their coadjoint
orbits as affine algebraic supervarieties. We find an algebraic
deformation quantization of them that can be related to the fuzzy
spaces of non commutative geometry.
\end{abstract}

\vfill\eject

\section{Introduction}

The use of   ``odd variables" in physics is very old and quite
natural. It is unavoidable in the description of theories
involving particles like the electron. Physical theories are given
by functionals on a certain kind of ``fields" or sections of a
bundle over the manifold of spacetime $\mathcal{M}$. The bundle
has fibers which are super vector spaces ($\Z_2$-graded vector
spaces), and that have a commutative superalgebra structure. In an
open set $U\subset \mathcal{M}$ of spacetime the set of fields can
be written as
$$\{\phi_i(x), \psi_\alpha(x)\}, \quad i=1,\dots p,\;
\alpha=1,\dots q,\quad x\in U,$$ with the superalgebra structure
given by the relations
\begin{eqnarray*}&\phi_i(x)\phi_j(x)-\phi_j(x)\phi_i(x)=0,\quad
\phi_i(x)\psi_\alpha(x)-\psi_\alpha(x)\phi_i(x)=0, \\
&\psi_\alpha(x)\psi_\beta(x)+\psi_\beta(x)\psi_\alpha(x)=0.\end{eqnarray*}

Free fields (i.e. fields satisfying the free equations of motion)
define unitary representations of the Poincar\'e group. It was
soon realized that in order to have a consistent description, the
half integer spin representations had to be described by odd
valued fields while the integer spin representations were
described by even valued fields. This is known as the
spin-statistics theorem.

Once the necessity of odd or Grassmann coordinates was established
on physical grounds it would have been natural to consider
superspaces or supermanifolds. Nevertheless, there was still an
important conceptual step to be taken, and the input came also
from physics. To have a supermanifold, all the coordinates, even
and odd, must be put on equal grounds. This means that changes of
coordinates mix together even and odd quantities in a consistent
way. From the physics point of view, this would mean that fields
with different spin and statistics could be mixed by certain
transformations. But all the groups of symmetries considered in
physics until then would preserve the even and odd subspaces, not
allowing such a mixing. They were all {\it even transformations}
or {\it super vector space morphisms}, since they would not change
the degree of the vector on which they act. The first example of
super transformations, that is, transformations that mix the odd
and even subspaces irreducibly, appeared in the seminal paper by
Gol'fand and Likhtman \cite{gl}. (Other examples of Lie
superalgebras had appeared in the mathematical literature before,
but not related to trasformation groups). In that paper,  an
extension of the Poincar\'e algebra by means of a set of odd
translations whose anticommutator is proportional to an ordinary
translation was constructed. To comply with the spin-statistics
theorem, the generators of the odd translations where in a spin
1/2 representation of the Poincar\'e group. This was the first
example in physics of a Lie superalgebra.

Natural as it seems from this historical perspective, the
introduction of Lie superalgebras (and Lie supergroups) has been
revolutionary in physics. The assumption that  {\it supersymmetry}
(that is, symmetry under a supergroup of transformations) is
realized in nature has profound implications, since it restricts
significantly the number of theories that can be considered. On
the other hand,   supersymmetric field theories have many
advantages in the quest for a unified theory of all interactions.
Accordingly, there is a huge number of works dealing with
supersymmetric field theories. Some pioneering works are the one
by  Wess and Zumino \cite{wz1} on the chiral multiplet, by Salam
and Strathdee \cite{ss} and Ferrara,  Wess and Zumino \cite{fwz}
on superspace, by  Wess and Zumino \cite{wz2} on supersymmetric
electrodynamics, by Ferrara and Zumino \cite{fz} and Salam and
Strathdee \cite{ss2} on supersymmetric Yang-Mills theories
\cite{fz}, and by Ferrara, Freedman and van Nieuwenhuizen
\cite{ffv} and Deser and Zumino \cite{dz} on supergravity. Further
information and references can be found in the collected reprint
volume \cite{fe}.

\bigskip

Since those first days the concept of supermanifold in mathematics
has evolved to a very precise formulation. Starting with the
superanalysis of Berezin \cite{be}, and  together with other works
by Kostant \cite{ko}, Leites \cite{le}, Manin \cite{ma1},
Bernstein, Deligne and Morgan \cite{dm} to mention some of the
most representative. Supermanifolds are seen roughly as ordinary
manifolds together with a sheaf of superalgebras. The odd
coordinates appear in the stalks of the sheaf. This approach
allows considerable freedom. Following Manin \cite{ma1}, one can
define different kinds of superspaces, supermanifolds or algebraic
supervarieties by choosing a base space with the appropriate
topology. The sheaves considered are superalgebra valued, so that
one can  generalize the concepts of complex and algebraic geometry
to this new setting.

In this paper
we will deal strictly with algebraic supervarieties, but the
concept of supermanifold has been treated extensively in the
literature mentioned above.

It is interesting to note that there is an alternative definition
of algebraic variety in terms of its functor of points.
Essentially, an algebraic variety can be defined as a certain
functor from some category of commutative algebras to the category
of sets. It is then very natural to substitute the category of
algebras by an appropriate category of commutative superalgebras
and to call this a  supervariety. The same can be done for
supergroups, super Lie algebras, coadjoint orbits of supergroups
and other ``super" objects. The elegance of this approach cannot
hide the many non trivial steps involved in the generalization. As
an example, it should be enough to remember the profound
differences between the classifications of semisimple Lie algebras
and semisimple Lie superalgebras \cite{be,ka,ri}.

The purpose of this paper is to study the coadjoint orbits of
certain  supergroups, to establish their structure as algebraic
supermanifolds and to obtain a quantum deformation of the
superalgebras that represent them.

Recently there has been a growing interest in the physics
literature on non commutative spaces (also called ``fuzzy
spaces"). The idea that spacetime may have non commuting
coordinates which then cannot be determined simultaneously  has
been proposed at different times for diverse motives
\cite{ma,do,co,cds}.  In particular, these spaces have been
considered as possible compactification manifolds of string theory
\cite{cds}.

In general, we can say that a ``fuzzy" space is an algebra of
operators on a Hilbert space obtained by some quantization
procedure. This means that it is possible to define a classical
limit for such algebra, which will be a commutative algebra with a
Poisson structure, that is, a phase space. As it is well known,
the coadjoint orbits of Lie groups are symplectic manifolds and
hamiltonian spaces with the Kirillov Poisson structure of the
corresponding Lie group. There is a wide literature on the
quantization of coadjoint orbits, most of the works based on the
Kirillov-Kostant orbit principle which associates to every orbit
(under some conditions) a unitary representation of the group. Let
us consider as an example the sphere $S^2$, a regular coadjoint
orbit of the group SU(2). In physical terms, the quantization of
$S^2$ is the quantization of the spin; it is perhaps the most
classical example of quantization, other than flat space
$\R^{2n}$. In the approach of geometric quantization (see Ref.
\cite{vo} for a review) the sphere must have half integer radius
$j$  and to this orbit it corresponds the unitary representation
of SU(2) of spin $j$ (which is finite dimensional). Let us take
the classical algebra of observables to be the polynomials on the
algebraic variety $S^2$. Given that the Hilbert space of the
quantum system is finite dimensional, the algebra of observables
is also finite dimensional. It is in fact of dimension $(2j+1)^2$
and isomorphic to the algebra of $(2j+1)\times (2j+1)$-matrices
\cite{vg,ma}. After a rescaling of the coordinates, we can take
the limit $j\rightarrow \infty$ maintaining the radius of the
sphere constant. In this way the algebra becomes infinite
dimensional and all polynomials are quantized. This procedure is
most appropriately described by Madore \cite{ma}, being the
algebras with finite $j$ approximations to the non commutative or
fuzzy sphere.

The  {\it deformation quantization} approach \cite{bffls}  is
inspired in the correspondence principle: to each classical
observable there must correspond a quantum observable (operator on
a Hilbert space). This {\it quantization map} will induce a non
commutative, associative {\it star product} on the algebra of
classical observables which will be expressed as a power series in
$\hbar$, having as a zero order term the ordinary, commutative
product and as first order term the Poisson bracket. One can ask
the reverse question, and explore the possible deformations of the
commutative product independently of any quantization map or
Hilbert space. It is somehow a semiclassical approach. Very often,
the star product is only a formal star product, in the sense that
the series in $\hbar$ does not converge. It however has the
advantage that many calculations and questions may be posed in a
simpler manner (see \cite{sw} for the string theory application of
the star product).

The deformation quantization of the sphere and other coadjoint
orbits has been treated in the literature \cite{ca,acg,cg,fl,fl2}.
The immediate question is if deformation theory can give some
information on representation theory and the Kirillov-Kostant
orbit principle \cite{fr,mo}.

In Ref. \cite{fl} a family of algebraic star products is defined
on regular orbits of semisimple groups. The star product algebra
is isomorphic to the quotient of the enveloping algebra by some
ideal. It is in fact possible to select this ideal in such a way
that it is in the kernel of some unitary representation. The image
of this algebra by the representation map will give a finite
dimensional algebra. In this way we can obtain the fuzzy sphere
and other fuzzy coadjoint orbits from the star product. The star
product is then seen as a structure underlying and unifying all
the finite dimensional fuzzy algebras. The classical limit and the
correspondence principle are seen naturally in this approach.

In this paper we extend the approach of Ref. \cite{fl} to some
coadjoint orbits of some semisimple supergroups. The extension to
the super category is not straightforward and this is reflected in
the fact that we have to restrict ourselves to the supergroups
$\rSL_{m|n}$ and $\rOSp_{m|n}$, and orbits of elements with
distinct eigenvalues (see Theorem \ref{orbittheorem}) to obtain
the extension. This is by no means an affirmation that the
procedure could not be in principle be applied to more general
supergroups and orbits, but a consequence of the technical
difficulties involved.

Other treatments of the ``fuzzy supershere" can be found in Refs.
\cite{bkr,gr,kl}. A comparison with these approaches as well as a
relation of our star product with representation theory of
supergroups will be discussed elsewhere.

The paper is organized as follows: Section
\ref{algebraicsupervarieties} is dedicated to the definition of
algebraic supervariety and its functor of points. The definitions
and results that we mention are somehow scattered in the
literature and we want to give a comprehensive account here
\cite{le2,ma}. Section \ref{as} is dedicated to algebraic
supergroups and their associated Lie superalgebras in terms of
their functors. The correspondence between the algebraic Lie
supergroup and its Lie superalgebra is treated in detail. Most of
the definitions and theorems extend easily from the classical
case, but we are not aware of any reference where this account has
been done explicitly for the algebraic case (instead, the
differential case is better known). We show proofs when we think
it can help to read the paper. We also illustrate the abstract
definitions with the examples of $\rGL_{m|n}$ and $\rSL_{m|n}$. In
Section \ref{cos} we retrieve the coadjoint orbits of the
supergroups mentioned above as certain representable functors, and
then as affine algebraic supervarieties. Finally, in Section
\ref{dco} we present a deformation of the superalgebra that
represents the functor associated to the coadjoint orbit.

\section{Algebraic supervarieties \label{algebraicsupervarieties}and
superschemes} In this section we want to give a definition of
algebraic supervarieties and superschemes. Our description is
self-contained and very much inspired in the approach  of Refs.
\cite{ma1,dm}. We assume some knowledge of basic algebraic
geometry and of super vector calculus. We are especially
interested in the description of a supervariety in terms of its
functor of points.

Let $k$ be a commutative ring. All algebras and superalgebras will be
intended to be over $k$ and are assumed to be commutative unless
otherwise stated. If $A$ is a  superalgebra we will denote by
$A_0$ the even part and by $A_1$ the odd part, so $A=A_0+A_1$. Let
$\cI^{\mathrm{odd}}_A$ be the ideal in $A$ generated by the odd
part. The quotient $A^\circ=A/\cI^{\mathrm{odd}}_A$ is an ordinary
 algebra. Notice that $A$ is both, an $A_0$-module and an
$A^\circ$-module.

If $\cO_X$ is a sheaf of  superalgebras over a topological
space $X$, then $\cO_{X,0}$ and $\cO_{X}^\circ$ are sheaves of
 algebras with $$\cO_{X,0}(U)=\cO_{X}(U)_0,\qquad
\cO_{X}^\circ(U)=\cO_{X}(U)^\circ\qquad U\subset_\mathrm{open}X.$$

We denote by $\arf$ the category of commutative $k$-algebras which
are reduced\footnote{A reduced algebra is an algebra that has no
nilpotent elements.}  and finitely generated (often called {\it
affine} algebras), and by $\sarf$ the category of commutative
$k$-superalgebras finitely generated, and such that, modulo the
ideal generated by the odd elements,  they are reduced.

We denote also by $\alg$ and $\sa$ the categories of commutative
$k$-algebras and commutative $k$-superalgebras respectively.

\subsection{Ringed superspaces}

\begin{definition}
A  superspace  $(X,\cO_X)$ is a topological space $X$ together
with a sheaf of superalgebras $\cO_X$, such that:

\noindent a. $(X, \cO_{X,0})$ is a noetherian scheme,

\noindent b. $\cO_X$ is a coherent sheaf of $\cO_{X,0}$-modules.
\end{definition}

\begin{definition} A morphism of superspaces $(X,\cO_X)$ and
$(Y,\cO_Y)$ is given by a pair $(f,\psi)$ where $f:X \rightarrow
Y$ is a continuous map, $\psi:\cO_Y \rightarrow f_*\cO_X$
is a map of sheaves of superalgebras on $Y$ and
$(f,\psi|_{\cO_Y})$ is a morphism of schemes.
\end{definition}
$f_*\cO_X$ denotes the push-forward of the sheaf $\cO_X$ under
$f$.

\begin{example}\label{afsu} \end{example}

Let $A$ be an object of $\sa$. We consider the topological space
$X_A:=\spec(A_0)=\spec(A^\circ)$ (they are isomorphic since the
algebras differ only by nilpotent elements) with the Zariski
topology.

On $X_A$ we have  the structural sheaf, $\cO_{A_0}$. The stalk of
the sheaf at the prime $\bp\in \spec(A_0)$ is the localization of
$A_0$ at $\bp$. As for any superalgebra, $A$ is a module over
$A_0$, so we have indeed a sheaf of $\mathcal{O}_{A_0}$-modules
over $X_A$ with stalk the localization of the $A_0$-module $A$
over each prime $\bp\in \spec(A_0)$. It is a sheaf of
superalgebras that we will denote by $\tilde A$. $(X_A,\tilde A)$
is a superspace. (For more details on the
construction of the sheaf $\tilde M$, for a generic $A_0$-module
$M$, see \cite{ha} II \S 5).  \hfill$\blacksquare$

\begin{definition}
An affine  algebraic supervariety is a superspace
isomorphic to $(X_A,\tilde A)$ for some superalgebra $A$
in $\sarf$.
\end{definition}
$A$ is often called the {\it coordinate superalgebra} of the
supervariety. The affine algebraic supervarieties form a category
denoted by $\svaraff$. Given an affine algebraic supervariety
$(V,\cO_V)$ (or just $V$ for short) we have an ordinary affine
algebraic variety associated to it, $(V,\cO_V^\circ)$. It is
called the {\it reduced variety} of $V$ and denoted also by
$V^\circ$.

\bigskip

 It is straightforward to generalize the  construction  to an
arbitrary commutative superalgebra $A$. We then have  the
following

\begin{definition} An affine superscheme is a superspace
 isomorphic to $(X_A,\tilde A)$ for some
superalgebra $A$ in $\sa$.\end{definition}

The affine superschemes over $k$ form a category denoted by
$\sschemeaff$. As in the classical case, general supervarieties
and superschemes are defined as superspaces that are locally affine
supervarieties and affine superschemes.

To any  superscheme $(S, \cO_S)$ one  can associate an ordinary
scheme $(S, \cO_S^\circ)$, called the {\it reduced scheme}.

\begin{example} Supervariety over the sphere $S^2$.\end{example}
(See Ref.\cite{be}, page 8). This is an example of a supervariety
explicitly given by the superalgebra representing the functor. We
will see how the sheaf is constructed in the closed points
(maximal ideals) of the topological space. Consider the free
commutative superalgebra generated by three even variables $x_1,
x_2, x_3$ and three odd variables $\xi_1, \xi_2, \xi_3$, $k[x_1,
x_2, x_3,\xi_1, \xi_2, \xi_3],$ and the ideal $$\cI=(x_1^2+ x_2^2+
x_3^2-1,\; x_1\cdot \xi_1+ x_2\cdot\xi_2+ x_3\cdot\xi_3).$$ The
superalgebra $k[V]=k[x_1, x_2, x_3,\xi_1, \xi_2, \xi_3]/\cI$
represents a supervariety whose reduced variety $V^\circ$ is the
sphere $S^2$. At each maximal ideal in $k[V]_0$,
$$
\bmax=(x_i-a_i,\, \xi_i\xi_j)\quad \mathrm{with}
\;\; i,j=1,2,3, \;\; a_i\in k\;\; \mathrm{and}\;\;
a_1^2+a_2^2+a_3^2=1,
$$
the local ring of $k[V]_0$ is the ring of
fractions
$$
(k[V]_0)_{\bmax}=\{\frac{f}{g}\;/\;
f,g\in k[V]_0, \; g\notin \bmax\}.$$ The stalk of the
structural sheaf at $\bmax$ is the localization of
$k[V]$ as a $k[V]_0$-module, that is
$$k[V]_{\bmax}=\{\frac{m}{g}\;/\; m\in k[V],\,
g\in k[V]_0, \; g\notin \bmax\}.$$ Notice that if
$a_1\neq 0$ (not all $a_i$ are zero simultaneously), then $x_1$ is
invertible in the localization and we have $$\xi_1=\frac
1{x_1}(x_2\xi_2+x_3\xi_3),$$ so $\{\xi_2,\xi_3\}$ generate
$k[V]_{\bmax}$ as an
$\cO_{k[V]_0}$-module.\hfill$\blacksquare$

\subsection{The functor of points}

 We recall first the definition of the {\it functor of points} in
 the classical (non super) case.

 Let $X$ be an affine variety. The representable functor
$h_X:\varaff^{\mathrm{opp}} \rightarrow \sets$ \footnote{By
 the  label opp  we denote the category with the direction of
 morphisms inverted. We could equally speak about contravariant functors.}
 from the category of affine varieties to the category of sets
$$h_X(Y)=\hom_{\mathrm{\varaff}}(Y,X)$$ is the functor of points
of $X$. An element of $\hom_{\mathrm{\varaff}}(Y,X)$ is an
$Y$-point of $X$.
 Given the  equivalence of categories
\begin{eqnarray*}F:&\arf&\longrightarrow\; \varaff\\& R &\longrightarrow \; \mathrm{Spec}(R),\end{eqnarray*}
we can equivalently define the functor of points of an affine
variety $X=\spec (R)$ as the  representable functor
$$h_R(T)=\hom_{\mathrm{\arf}}(R,T), \qquad T\in \arf.$$ This
leads to an alternative definition of affine variety as a
representable functor between the categories $\arf$ and
$\sets$.

\bigskip

These definitions and observations can be extended immediately to
the super case.

\begin{definition}
The functor of points  of a supervariety $X$ is a representable
functor
$$
h_X:\svaraff^{\mathrm{opp}}\rightarrow \sets.
$$
\end{definition}

\begin{observation} \end{observation}
The construction of Example \ref{afsu} defines an equivalence of
categories between $\svaraff$ and $\sarf$ and between
$\sschemeaff$ and $\salg$. So the functor of points can be
equivalently be given by a representable functor
$h_A:\sarf\rightarrow \sets$. For example, as we will see in the
next section,  algebraic supergroups can easily be defined as
certain representable functors.

 \bigskip

In the general case, for a general  superscheme $X$,  its functor
of points can be defined either as a representable functor from
$\sschemes^\mathrm{opp} \rightarrow \sets$ or  as the functor
$$
h_X:\sa^{\mathrm{opp}} \lra \sets, \qquad
h_X(A)=\hom_{\sschemes}(X_A,X)
$$
where $X_A$ is the affine superscheme associated to the
superalgebra $A$ (see Example \ref{afsu}). In fact, as it happens
in the non super case (see Ref. \cite{eh} pg 253), the functor of
points of a superscheme is determined by looking at its
restriction to  affine superschemes. Since we are mainly
interested in affine supervarieties and affine superschemes we
will not pursue further this subject.

\section{ \label{as}Algebraic supergroups and their Lie superalgebras}
In this section we generalize the basic notions of algebraic
 groups to the super case.
\subsection{\label{ssf}Supergroups and supergroup functors}

\begin{definition}
An affine algebraic supergroup is a
representable, group valued functor
$$
G:\sarf\rightarrow \sets.
$$
\end{definition}

The superalgebra $k[G]$ representing a supergroup has the
additional structure of a commutative super Hopf algebra. The
coproduct
\begin{eqnarray*}\Delta:k[G]&\longrightarrow &k[G]\otimes k[G]\\
f&\longrightarrow &\Delta f\end{eqnarray*} is  such that for any
affine superalgebra $A$ and two morphisms ($A$-points) $x,y\in
G(A)=\hom_{k-\mathrm{superalg}}(k[G],A)$ the following relation
holds: $$m_A\bigr(x\otimes y(\Delta f)\bigl)=x\cdot y(f) \qquad
\forall f\in k[G],$$ where $m_A$ denotes the multiplication in $A$
and $``\cdot"$ denotes the multiplication in the group $G(A)$. The
tensor product is understood in the tensor category of super
vector spaces \cite{dm}. The counit is given by
\begin{eqnarray*}\E:&k[G]\longrightarrow &k\\
&f\longrightarrow &e(f)\end{eqnarray*} where $e$ is the identity
in $G(k)$. The antipode is defined as
\begin{eqnarray*}S:&k[G]\longrightarrow &k[G]\\
&f\longrightarrow &S(f)\end{eqnarray*} with
$x(f)=x^{-1}(S(f))\quad \forall x\in G(A),$ and where $x^{-1}$
denotes the group inverse of $x$ in $G(A)$.

\begin{example} \label{supergroups}
The supergroups $\rGL_{m|n}$ and $\rSL_{m|n}$.\end{example}
Let $A$ be a commutative superalgebra. We denote by $A^{m|n}$ the
free module over $A$ generated by $m$ even generators and $n$ odd
generators. The endomorphisms of this super vector space (linear
maps that preserve the grading) are given by matrices (we use the
conventions of Ref. \cite{dm})
\begin{equation}\begin{pmatrix}p_{m\times m}&q_{m\times n }\\r_{n\times
m}&s_{n\times n}
\end{pmatrix}\label{morphisms}\end{equation}
where $p$ and $s$ have even entries  and $q$ and $r$ have odd
entries in $A$. We denote the set of these matrices as
$\fgl_{m|n}(A)$. We can define a functor
\begin{eqnarray*}\fgl_{m|n}&:\sarf&
\longrightarrow \sets\\ &A&\longrightarrow
\fgl_{m|n}(A).\end{eqnarray*} $\fgl_{m|n}$ is a representable
functor. It is represented by the superalgebra
$$k[\fgl_{m|n}]:=k[x_{ij}, y_{\alpha \beta},
\xi_{i\beta},\gamma_{\alpha j}],\qquad i,j=1,\dots m,\;\;
\alpha,\beta=1,\dots n$$ where $x$ and $y$ are even generators and
$\xi$ and $\gamma$ are odd generators. In fact, writing the
generators in matrix form, $$\begin{pmatrix}
\{x_{ij}\}&\{\xi_{i\beta}\}\\\{\gamma_{\alpha j}\}& \{y_{\alpha
\beta}\}&\end{pmatrix},
$$
any matrix as in (\ref{morphisms}) assigns to a generator of
$k[\fgl_{m|n}]$ an element of $A$, and the assignment has the
right parity. Hence, it defines a superalgebra morphism
$k[\fgl_{m|n}]\rightarrow A$.

\medskip

 $\rGL_{m|n}(A)$ is defined as the set of all morphisms $g:A^{m|n}\rightarrow
A^{m|n}$ which are invertible.
In terms of the matrix ({\ref{morphisms}), this means that the
{\it Berezinian} \cite{be} or {\it superdeterminant} $$\Ber(g)=
\det(p-qs^{-1}r)\det(s^{-1})$$ is invertible in $A$. A
necessary and sufficient condition for $g$ to be invertible is that
$p$ and $s$ are invertible. The group valued functor
\begin{eqnarray*}
\rGL_{m|n}&:\sarf &
\longrightarrow \sets\\ &A&\longrightarrow
\rGL_{m|n}(A).
\end{eqnarray*}
is an affine supergroup represented by the algebra \cite{fi}
\begin{eqnarray*}&
k[\rGL_{m|n}]:=k[x_{ij}, y_{\alpha \beta},
\xi_{i\beta},\gamma_{\alpha
j},z,w]/\bigr((w\det(x)-1,z\det(y)-1\bigl),
\\
&i,j=1,\dots m,\;\;
\alpha,\beta=1,\dots n.
\end{eqnarray*}
Requiring that the
berezinian is equal to 1 gives the supergroup $\rSL_{m|n}$,
represented by
$$
k[\rGL_{m|n}]/\bigr(\det(x-\xi y^{-1} \gamma)z-1\bigl).
$$
where $y^{-1}$ is the matrix of indeterminates, inverse
of the matrix $y$, whose determinant is invertible and has inverse $z$.

In Ref.\cite{fi} the Hopf superalgebra
structure of this affine supergroup was explicitly computed. \hfill$\blacksquare$

\bigskip

 In the
classical case, the concept of {\it group functor} is a
generalization of the concept of algebraic group. This is treated
extensively in \cite{dg} II, \S 1.
It will be useful to
introduce this notion for the super case, which can be done
easily with suitable changes.

\begin{definition}
Let $G$ be a functor from $\sa$ to $\sets$. We say that $G$ is a
 supergroup functor if:

1. There exists a natural transformation called the {\it
composition law}
$$
m:G \times G \lra G
$$
satisfying the associativity condition: $m \cdot (m \times \id)=m
\cdot (\id \times m)$.

2. There exists a natural transformation, the {\it unit section},
$u:e_k \lra G$, where $e_k:\sa \lra \sets$, $e_k(A)=1_A$,
satisfying the commutative diagram:
$$
\begin{array}{ccccc}
G \times e_k & \stackrel{id \times u}\lra &
G \times G & \stackrel{u \times id}\lra & e_k \times G \\
             & _{}\searrow & _{m}\Big\downarrow & \swarrow  _{}& \\
             &       & G & & \\
\end{array}
$$

3. There exists a natural transformation $\sigma:G \lra G$
satisfying:
$$
\begin{CD}
G @> (id, \sigma)>> G \times G \\
@V{}VV @VVmV \\
  e_k @>u>> G  \\
\end{CD}
$$
\end{definition}

\medskip

A morphism of supergroup functors is defined as a natural
transformation preserving the composition law.

Because of their representability property, affine algebraic
supergroups are supergroup functor.

\subsection{Lie superalgebras}

Let $\cO_k:\sarf \lra \sets$ be the functor represented by $k[x]$.
 $k[x]$ corresponds in fact  to an ordinary algebraic variety, $\spec(k[x])$, the affine
line. For a superalgebra $A$ we have that $\cO_k(A)=A_0$.

\begin{definition}\label{superalgebra}
A Lie superalgebra is a representable, group valued functor
$$
\fg:\sarf \lra \sets
$$
with the following properties:

1. $\fg$ has the structure of $\cO_k$-module, that is, there is a
natural transformation $\cO_k \times \fg \lra \fg$. For each
superalgebra $A$ we have an  $A_0$-module structure on $\fg(A)$.

2. There is a natural transformation $[\;,\;]:\fg \times \fg \lra
\fg$ which is $\cO_k$-linear and that satisfies commutative
diagrams corresponding to the antisymmetric property and the
Jacobi identity. For each superalgebra $A$, $[\;,\;]$ defines a
Lie algebra structure on $\fg(A)$, hence the functor $\fg$ is Lie
algebra valued.
\end{definition}

For any algebraic supergroup there is a Lie superalgebra which is
naturally associated. It is our purpose to construct it
explicitly. Again, the construction is a generalization (not
completely straightforward) of what happens in the non super case.
Our treatment is very similar to the one in Ref.\cite{dm} II, 4,
no 1.

\medskip

Let $A$ be a commutative superalgebra and let
$A(\e)=_{\defi}A[\e]/(\e^2)$ be the algebra of dual numbers ($\e$
here is taken as an \emph{even} indeterminate). We have that
$A[\e]=A \oplus \e A$ and there are two homorphisms: $i: A \ra
A(\e)$ defined by $i(1)=1$ and  $p: A(\e) \ra A$ defined by
$p(1)=1$, $p(\e)=0$.

\begin{definition}
Let $G$ be a supergroup functor. Consider the homomorphism
$G(p):G(A(\e)) \longrightarrow G(A)$. For each $G$ there is a
supergroup functor, $\Lie(G)$, defined as
$$
\Lie(G)(A)=_{\defi}\ker(G(p)).
$$
\end{definition}
 We will see that when $G$ is an affine algebraic supergroup (so
 it is a representable group functor) $\Lie(G)$ is indeed a Lie superalgebra.

It is instructive to see first an example.

\begin{example} $\Lie(\rGL_{m|n})$, $\Lie(\rSL_{m|n})$ \end{example}
We want to determine the functor $\Lie(\rGL_{m|n})$.
Consider the map:
\begin{equation*}
\begin{array}{cccc}
\rGL_{m|n}(p):  & \rGL_{m|n}(A(\e)) &  \lra  & \rGL_{m|n}(A)
\\[1ex]
& {\begin{pmatrix} p+\e p' & q+\e q' \\ r+\e r' & s+ \e s' \end{pmatrix}}
& \mapsto &
{\begin{pmatrix} p & q \\ r & s \end{pmatrix}}
\end{array}
\end{equation*}
with $p,p',s,s'$ having entries in $A_0 $ and $q,q',r,r'$ having
entries in $A_1$.  $p$ and  $s$ are invertible matrices. One can
see immediately that
$$
\Lie(\rGL_{m|n})(A)=\ker(\rGL_{m|n}(p))= \bigl\{\begin{pmatrix}
1+\e p' & \e q' \\ \e r' & 1+ \e s' \end{pmatrix} \bigr\}.
$$
The functor $\Lie(\rGL_{m|n})$ is clearly group valued and can be
identified with the (additive) group functor $\fgl_{m|n}$ (see
Example \ref{supergroups}).

For $\Lie(\rSL_{m|n})$ one gets the extra condition $
\Ber=\det(1+\e p')\det(1-\e s')=1 $, which implies  the zero
supertrace condition, that is, $tr(p')-tr(s')=0$. The functor
$\Lie(\rSL_{m|n})$ is then identified with the (additive)
group functor
$$
\fsl_{m|n}(A)=\{x \in \fgl_{m|n}(A)\;/\; \str(x)=0\}.
$$

The functors $\fgl_{m|n}$ and $\fsl_{m|n}$ are representable and
Lie algebra valued.  We have already seen that $\fgl_{m|n}$ is
representable (see Example \ref{supergroups}). $\fsl_{m|n}$ is
represented by the superalgebra:
$$
k[\fsl_{m|n}]=k[\rM_{m|n}]/\str(m)
$$
where $m$ is the matrix of indeterminates generating the algebra
$k[\fgl_{m|n}]$.\hfill$\blacksquare$

\medskip

We want to show that when $G$ is an affine algebraic group then
$\Lie(G)$ is a superalgebra. We will show first that the functor
satisfies both properties in Definition \ref{superalgebra}. Then
we will show that it is representable.  We start by seeing  that
$\Lie(G)$ has a structure of $\cO_k$-module.

 Let $u_a:A(\e) \lra A(\e)$ be the endomorphism,
$u_a(1)=1$, $u_a(\e)=a\e$, for $a \in A_0$. $\Lie(G)$ admits a
$\cO_k$-module structure, i.e. there is a natural transformation
$$
\begin{CD}
\cO_k \times \Lie(G) @>>> \Lie(G)
\end{CD}
$$such that for any superalgebra $A$, $a\in \cO_k(A)$, $x\in \Lie(G)(A)$, $$
(a,x) \mapsto ax=\Lie(G)(u_a)x.$$ Notice that for subgroups of
$\rGL_{m|n}(A)$, $ax$ corresponds to the multiplication of the
matrix $x$ by the
 even scalar $a$.

 Let  us consider the group of linear automorphisms of
 $\Lie(G)(A)$. Because of the natural $\cO_k$-module structure of
 $\Lie(G)$ we have a group functor
 $$\rGL(\Lie(G)):\sarf \lra \sets.$$
 Alternatively one can   also denote
 $\rGL(\Lie(G))=\mathrm{Aut}(\Lie(G))$.

We now want to introduce the Lie algebra structure (property 2 in
Definition \ref{superalgebra}). We will do it through the adjoint
actions, seen as natural transformations.

\begin{definition}
Let $G$ be a supergroup functor. The adjoint action of $G$ on
$\Lie(G)$ is defined as the natural transformation
$$
\begin{CD}
G @ > \Ad>> \rGL(\Lie(G)) \end{CD}$$ which for any superalgebra
$A$  and $g \in G(A), \;x \in \Lie(G)(A)$
$$\Ad_A(g)x=_{\defi}G(i)(g)xG(i)(g)^{-1}
$$
\end{definition}
 $\Ad_A(g)x\in G(A(\e))$ but  since $G(p)$ is a group homomorphism,
$\Ad_A(g)x\in\Lie(G)(A)$.

 $\Ad$  is a morphism in the category of group functors. So we can
 make the following definition:
\begin{definition} Let $G$ be a supergroup functor. The adjoint
action of $\Lie(G)$ on $\Lie(G)$ is defined as
$$
\ad=_{\defi}\Lie(\Ad):\Lie(G) \lra
\Lie(\rGL(\Lie(G)))=\End(\Lie(G))
$$\end{definition}
Finally, the  bracket on $\Lie(G)(A)$ is defined as
$$
[x,y]=_{\defi}\ad(x)y, \qquad x,y \in \Lie(G)(A).
$$
The arguments in Ref.\cite{dg} II \S 4 4.2, 4.3 apply with small
changes to our case and prove that $[\;,\;]$ is a Lie bracket.
\medskip





\begin{example}$\rGL_{m|n}$\end{example}
We want to see that in the case of $\rGL_{m|n}$ the Lie bracket
$[\;,\;]$ is the commutator. We have
$$
\begin{array}{c}
\begin{array}{cccc}
\Ad_A:& \rGL(A) & \lra & \rGL(\Lie(\rGL_{m|n}))(A)=\rGL(\rM_{m|n}(A)) \\
& g & \mapsto & \Ad(A)(g),
\end{array}
\\ \\
\Ad(g)x=gxg^{-1}, \qquad x\in \rM_{m|n}(A).
\end{array}
$$
By definition we have: $\Lie(\rGL(\fgl_{m|n}))(A)=\{ 1+\e \alpha
\;/\; \alpha \in \rGL(\fgl_{m|n})(A) \}$
So we have, for $a,b \in \fgl_{m|n}(A) \cong \Lie(\rGL_{m|n})(A)=$
$\{1+\e a\;|\;a \in \fgl_{m|n}(A)\}$:
$$
\ad(1+\e a)b=(1+\e a)b(1-\e a)=b+(ab-ba)\e=b+\e [a,b].
$$
Hence $\ad(1+\e a)=\id+\e \alpha(a)$, with $\alpha(a)= [a,\;]$ as
we wanted to show.

It is also clear that the same computation will hold for any
closed subgroup of $\rGL_{m|n}$. \hfill$\blacksquare$

\medskip

Finally we have to address the issue of the  the representability
of the functor $\Lie(G)$. From the classical (non super) case, we
know that it is not, in general, representable (see Ref.\cite{dg}
II, 4, 4.8). The most useful example is  when $G$ is a linear
group, that is, a closed subgroup of $\rGL_n$. Then $\Lie(G)$ is
representable, and we will see that the same is true for the super
case.

More generally, for any affine algebraic supergroup $G$, the
associated functor $\Lie(G)$  can be shown to be representable.
Theorem \ref{lierep} characterizes $\Lie(G)$ geometrically as the
tangent space at the identity. In the classical case, this result
is a particular case of a more general one involving schemes,
whose proof can be found in Ref. \cite{dg}. Since we are only
concerned with
 affine algebraic supergroups, all we need is an extension  to the
 super category of the (simpler) proof for affine algebraic groups.
 We find then convenient to write explicitly the proof of the following theorem, although
  the super extension in this particular case presented no difficulty.

\begin{theorem} \label{lierep} Let $k$ be a field and $G$ be an affine algebraic group, with
 $k[G]$ its coordinate superalgebra.
As in Section \ref{ssf}, let $\E$ denote the counit in the
superHopf algebra $k[G]$. We denote by  $\m_\ep=\ker(\ep)$ and by
$\omega$ be the super vector space $\m_\ep/\m_\ep^2$. Then,
$\Lie(G)$ is a representable functor and it is represented by
$k[\omega]$ where
$$k[\omega]=k[x_1, \dots x_p, \xi_1, \dots \xi_q]$$ with $x_1 \dots
x_p$, and  $\xi_1 \dots \xi_q$ being even and odd indeterminates
respectively, and $p|q$ is the superdimension of $\omega$.
\end{theorem}

 {\it Proof}. We have to prove that
$\Lie(G)(A)=\hom_{k-\mathrm{superalg}}(k[\omega],A) $. It is
immediate to verify that
$$\hom_{k-\mathrm{superalg}}(k[\omega],A)\cong\hom_{k-\mathrm{supermod}}(\omega,A)\cong
(\omega^* \otimes A)_0,$$ hence it is enough to  show that
$\Lie(G)(A)\cong\hom_{k-\mathrm{supermod}}(\omega,A)$.

We will define a map $\rho: \hom_{k-\mathrm{supermod}}(\omega,A)
\rightarrow \Lie(G)(A)$ and then show that it is a bijection.
 Let $d\in
\hom_{k-\mathrm{supermod}}(\omega,A)$, $d:\omega \lra A$. We first
consider the following maps:
$$
\begin{array}{cccc}
\phi:& k[G] & \lra & k\oplus \omega=k \oplus \m_\ep/\m_\ep^2 \\
& f & \mapsto & (\ep(f), f-\ep(f)+ \m_\ep^2),\\
\end{array}
$$
and for each $d$
$$
\begin{array}{cccc}
d':& k \oplus \omega & \lra & A(\e) \\
& (s,t) & \mapsto & s+d(t)\e .\\
\end{array}
$$
We define $\rho(d)$ as the composition $\rho(d)=d'\circ\phi$. Then
we have that
$$G(p)(\rho(d))(f)=\ep(f),$$
so $\rho(d)\in \ker G(p)=\Lie(G)(A)$.

We want now to give the inverse $z:\Lie(G)(A) \rightarrow
\hom_{k-\mathrm{supermod}}(\omega,A)$. Assume $\psi \in
\Lie(G)(A)$. We can write: $\psi(f)=\ep(f)+\e f'$. Now consider
$$
\begin{array}{cccc}
\psi|_{\m_e}:& \m_e & \lra & A(\e)\\
& f & \mapsto & \e f'
\end{array}
$$
Observe that $\psi|_{\m_e^2}=0$, hence going to the quotient we
have a supermodule map $\tilde \psi:\omega \lra A(\e)$,
$\tilde\psi(f)=f'\e$. Now define $z(\psi)(f)=f'$. This is the
inverse of $a$. The fact $z \cdot \rho=id$ is straightforward. For
$\rho \cdot z$ observe that given $\psi:k[G] \lra A(\e)$, this can
always be written as: $\psi(g)=\ep(g)+g''\e$. Observe that
$g-\ep(g) \in \m_e$ hence $\psi|_{\m_e}(g-\ep(g))=g''\e$, now the
result follows easily. \hfill$\blacksquare$.
\bigskip
\begin{corollary} Let $k$ be a field and
$G$ be an affine algebraic supergroup. Then $\Lie(G)$ is a
Lie superalgebra
\end{corollary}\hfill$\blacksquare$

\bigskip

We want to remark that $\omega=\m_e/\m_e^2$ can be regarded as the
dual of the tangent space at the identity of the supergroup $G$.
Theorem \ref{lierep} hence states that such tangent space is the
same as $\Lie(G)$, as it happens in the non super case.

\begin{observation} Lie superalgebras as super vector spaces with a graded bracket.
 \end{observation}

Lie superalgebras were first introduced in physics \cite{gl} with
a different definition. A Lie superalgebra is a  super
($Z_2$-graded) vector space $\fg=\fg_0+\fg_1$ with a bilinear,
graded operation
\begin{eqnarray*}[\;,\;]:&\fg\otimes \fg\longrightarrow &\fg\\
&X\otimes Y\longrightarrow &[X,Y]\end{eqnarray*} such that

\smallskip

\noindent 1. $[X,Y]=-(-1)^{p_Xp_Y}[Y,X]$

 \noindent 2. $[X,[Y,Z]]+
 (-1)^{p_Xp_Y+p_Xp_Z}[Y,[Z,X]]+(-1)^{p_Xp_Z+p_Yp_Z}[Z,[X,Y]]=0$

\noindent where $X,Y,Z$ are homogeneous elements of $\fg$ with
parities $p_X,p_Y,p_Z$. This definition is used mostly when Lie
superalgebras are treated independently of Lie supergroups
\cite{be,ka}.

This definition of Lie superalgebra can be shown to be equivalent
to Definition \ref{superalgebra} in functorial terms. This is
proven in Ref. \cite{dm}, Corollary 1.7.3 pg 57, using the {\it
even rules} principle. We will  show how it works for the specific
example $\fgl_{m|n}$ (Example \ref{exsuper}).

 We have then that in the super vector space
$\omega^*$ there is a graded bracket $[\;,\;]$ and a Lie
superalgebra structure in the sense mentioned above.

\begin{example}\label{exsuper}\end{example}
Let us consider the set of $(m+n)\times (m+n)$ matrices with
arbitrary entries in $A$ denoted as $\overline\fgl_{m|n}(A)$. We
also denote
 \begin{eqnarray*}I=1,\dots m+n, \quad && I=i\;\;\mathrm{for}\;\;I=1,
 \dots m\\&& I=n+\alpha\;\;\mathrm{for} \;\; I=n+1,\dots
 m+n.\end{eqnarray*}
 Let $\{E_{IJ}\}$ be the standard basis
of matrices with 1 in the place $IJ$ and 0 everywhere else. An
element $X\in \overline \fgl_{m|n}(A)$ can be written in terms of
the standard basis
$$X=X_{IJ}E_{IJ}=p_{ij}E_{ij}+q_{i\beta}E_{i\beta}+r_{\alpha
j}E_{\alpha j} + s_{\alpha \beta}E_{\alpha \beta},$$ where sum
over repeated indices is understood and the parities of $p,q,r$
and $s$ are arbitrary.

 We assign even degree to  $E_{ij}$ and $E_{\alpha\beta}$ (block diagonal matrices)
  and odd degree
 to $E_{i\beta}$ and $E_{\alpha j}$ (block off diagonal matrices). This corresponds
 to even and odd
 linear maps. With this assignments $\overline \fgl_{m|n}(A)$ is a non
 commutative,
 associative
superalgebra, and $\fgl_{m|n}(A)$ is its even part. It corresponds
to the even linear maps or super vector space morphisms.

   We can give it a super Lie algebra structure with
the ordinary commutator of matrices among even elements or among
an even and an odd one and the ordinary anticommutator of matrices
among odd elements. Then, the even part of this Lie superalgebra
is the Lie algebra $\fgl_{m|n}(A)$.

In general, giving a representable, Lie algebra valued functor is
equivalent to give a super Lie algebra through the
(anti)commutation rules of the generators.\hfill$\blacksquare$

\section{Coadjoint orbits of supergroups.\label{cos}}

Let $G \subset \rGL(m|n)$ be an algebraic Lie supergroup and
$\fg=\Lie(G)$ the associated Lie superalgebra. Let $\fg^*$ be the
functor $\fg^*(A)=\fg(A)^*$.

We want  to define a coadjoint orbit of the supergroup $G$. Let
$X_0$ be a geometric point of $\fg^*$, that is $X_0\in \fg^*(k)$.
$\fg(k)$ is an ordinary Lie algebra over $k$, and $\fg(k)\subset
\fg(A)$ for any $A$ through the unit map of $A$. We consider the
following functor,
\begin{eqnarray}
\cC_{X_0}:&\sa &\longrightarrow \sets\nonumber\\
&A&\longrightarrow \{\Ad_g^*X_0, \forall g\in
G(A)\}=G(A)/H(A)\label{coadfunctsup},
\end{eqnarray}
where $H(A)$
is the stability group of $X_0$.  Notice that it is necessary to
choose a geometric point in order to have a functor.

The functor $\cC_{X_0}$ is not, in general, representable. The
problem arises already at the classical level. We can see it with
an example. Let us consider the algebraic group $\rSL_n$ over the
complex numbers and its Lie algebra $\fsl_n$. Their functors of
points are represented respectively by
$$
\begin{matrix}
\C[\rSL_n]=\C[x_{ij}]/\bigr(\det(x)-1\bigl),\\ \\
\C[\fsl_n]=\C[x_{ij}]/\bigr(\tr(x)\bigl)
\end{matrix} \qquad
i,j=1,\dots n .
$$
Let $X_0=\diag(l_1 \dots l_{n}) \in \fsl_n(\C)$, $l_i \neq l_j$,
for $i \neq j$. The coadjoint orbit of $X_0$ is an algebraic
variety represented by
\begin{equation}
\C[\cC_{X_0}]=\C[x_{ij}]/(p_1-c_1,\cdots p_l-c_l),\label{repco}
\end{equation}
where $p_i=\tr(X^i)$ and $c_i=p_i(X_0)$.

An $A$-point of $\cC_{X_0}$ is a morphism
$\C[\cC_{X_0}]\rightarrow A$ and it is given by a matrix $(a_{ij})
\in \fgl_{n}(A)$, such that $p_k(a_{ij})=c_k$, $k=1 \dots m+n$.

If the functor $\cC_{X_0}^r$, defined as $\cC_{X_0}$ but
restricted to the category of commutative algebras were
representable, any such matrix would be of the form
$X=gX_0g^{-1},$ that is, conjugate to a diagonal one through an
element $g\in \rSL_n(A)$. But this is not necessarily true in an
arbitrary algebra $A$.

\bigskip

In the classical case the functor of points of the coadjoint orbit
is obtained as the {\it sheafification} (see Definition
\ref{sheafification}) of $\cC_{X_0}^r$ (\cite{dg} pg 341). We will
see that the same is true for the super case.

We will start by considering a functor $$ \cF:\sa
\longrightarrow \sets.
$$
This is a generalization of the concept of $\Z$-functor defined in
\cite{dg} pg 9.

 $\cC_{X_0}$ is an example of  such a functor.
Due to the equivalence of categories between $\sa$ and
$\sschemeaff$ we can give equivalently the functor:
$$
\cF:\sschemeaff^\mathrm{opp} \longrightarrow \sets.
$$
(Abusing the notation, we use the same letter for both functors).

For each affine superscheme $X$ the restriction $\cF_X$ of $\cF$
to $\Top(X)$ defines a presheaf on $X$ ($\Top(X)$ denotes the
category of open sets in $X$). A functor $ \cF:\sa \longrightarrow
\sets$ is said to be {\it local} if $\cF_X$ is a sheaf for each $X
\in \sschemes$.

Equivalently,  $ \cF$ is {\it local} if for any $A \in \sarf$ and
elements $f_1 \dots f_s \in A$ such that $(f_1 \dots f_s)=(1)$ we
have the exact sequence:
$$
F(A) \;\stackrel{F(\alpha_i)}\rightarrow\;
F(A_{f_i})\;\stackrel{F(\alpha_{ij})}\rightrightarrows\;
F(A_{f_if_j})
$$

In particular, any representable functor is local.

\begin{definition} \label{sheafification}Let us consider a functor
 $ \cF:\sa \longrightarrow
\sets$. We will denote by $\tilde \cF$ the unique local functor
such that for $f_1 \dots f_s \in A$, $(f_1 \dots f_s)=(1)$,
$\tilde \cF(A_{f_i})=\cF(A_{f_i})$. $\tilde \cF$ is called the
sheafification of $\cF$.\end{definition}

It is perhaps more natural to introduce the sheafification in
terms of the presheaves $\cF_X$. To any presheaf one can associate
a sheaf that is the ``closest" sheaf to the given presheaf; it is
called its sheafification. Then  $\tilde \cF$ can be equivalently
defined as the local functor obtained by doing the sheafification
of the presheaves $\cF_X$ for all $X \in \sschemes$.

\begin{observation}. \label{obsfunctor}
If $\cP$ is a subfunctor of a local functor $\cF$ and for all $R
\in \sa$ there exist $f_1 \dots f_r \in R^o$ such that $(f_1 \dots
f_r)=(1)$ and $\cP(R_{f_i})=\cF(R_{f_i})$ then $\tilde \cP$=
$\cF$. ($R$ is viewed as an $R_0$ module).
\end{observation}

\bigskip

Our strategy to determine the functor of
the coadjoint orbits of supergroups will be to find a
representable  (hence local) functor and to prove that it is the
sheafification of (\ref{coadfunctsup}).

We need first some notions about the invariant polynomials of a
superalgebra.

\subsection{Invariant polynomials}

For the rest of the section we take $k=\C$. The Cartan-Killing
form of a super Lie algebra $\fg$ is a natural trasformation:
$B:\fg \times \fg \lra \fg$,
$$
B_A(X,Y)=\str(\ad_X\ad_Y)
\qquad X, Y \in \fg(A).
$$
$B_A$ is an invariant, supersymmetric bilinear form. In the
following, we will consider simple Lie superalgebras whose
Cartan-Killing form is non degenerate, namely $\fsl_{m|n}$ with
$m\neq n$, $\fosp_{m|n}$ with $\frac{m}{2}-\frac{n}{2}\neq 1$
($m,n$ even), $\ff_4$ and $\fg_3$ \cite{ka}.

There is a natural isomorphism between the functors $\fg$ and
$\fg^*$ such that to each object $A$ in $\sarf$ assigns a morphism
\begin{eqnarray*}
\varphi_A:&\fg(A)\longrightarrow&\fg^*(A)\\
&X\longrightarrow& \varphi_A(X)\;\; \hbox{with}\;\;
\varphi_A(X)(Y)=B_A(X,Y).
\end{eqnarray*}
This isomorphism gives also an isomorphism between the
superalgebras representing both functors, $\C[\fg]\simeq\C[\fg^*]$
and intertwines  the adjoint and coadjoint representation, so the
adjoint orbits are the same as the coadjoint orbits. From now on
we will use the algebra $\C[\fg]$.

Let $G$ be an affine algebraic supergroup,
 subfunctor of $\rGL(m|n)$,
with super Lie algebra $\fg$. We say that $p \in \C[\fg]$ is an
{\it invariant polynomial} if for any $A$-point
$x:\C[\fg]\rightarrow A$ of $\fg$ and $g:\C[G]\rightarrow A$ of
$G$ we have that
\begin{equation}
x(p)=\Ad_A(g)x(p). \label{invpol}
\end{equation}
The invariant polynomials are a subalgebra of $\C[\fg]$. Contrary
to what happens in the classical case, this algebra may be not
finitely generated \cite{be,ksls}. This is the case for the
algebra of invariant polynomials on $\fgl_{m|n}$. The generators
can be taken to be the supertraces of arbitrary order,
$\str(X^k)$, which are independent. The invariant polynomials in
the reduced Lie algebra, $\fgl_m\times\fgl_n$ are generated by
traces $\tr(X^k)$ with $k=1,\dots m+n$, since higher order traces
can be expressed in terms of the first $m+n$ ones.

\subsection{The coadjoint orbits of a supergroup as algebraic supervarieties}

For a regular, semisimple element $X_0$ of $\fg(\C)$, its orbit
under the adjoint action of the group $G(\C)$ is an algebraic
variety defined by the values of the homogeneous Chevalley
polynomials $p_1,\dots p_l$, where $l$ is the rank of the group.
We will see that the supersymmetric extensions of these
polynomials define the adjoint orbit of the supergroup. More
specifically, we have the following

\begin{theorem} \label{orbittheorem}
Let $\hat p_1,\dots \hat p_l$ be polynomials on a simple Lie
superalgebra  of the type $\fsl_{m|n}$ or $\fosp_{m|n}$ with the
following properties:

\smallskip

\noindent 1. They are invariant polynomials under the adjoint
action (\ref{invpol}).

\smallskip

\noindent 2. Let $ p_1,\dots p_l$ be the  projections of
$\hat p_1,\dots \hat p_l$ onto
the reduced algebra $\C[\fg^\circ]=\C[\fg]/\cI^{\mathrm odd}$. The
ideal of the orbit of a regular, semisimple element $X_0\in
\fg(\C)$ with distinct eigenvalues, is
$\cJ=(p_1-c_1,\dots p_l-c_l)$, $c_i=p_i(X_0)$. So the
orbit of $X_0$ is an algebraic
variety whose functor is represented by $$\C[\fg^\circ]/\cJ.$$

\smallskip
\noindent Then
the sheafification of the functor $\cC_{X_0}$ (\ref{coadfunctsup})
is representable and is represented by
\begin{equation}A_{X_0}=\C[\fg]/\cI.\label{representable}\end{equation}
with $\cI=(\hat p_1-c_1,\dots \hat p_l-c_l)$.
\end{theorem}

{\it Proof}. Let us denote by $F:\sarf\rightarrow\sets$ the
functor represented by $A_{X_0}$ in (\ref{representable}). It is clear
that $\cC_{X_0}$ is a subfunctor of $F$, since for any
superalgebra $R$, an element of  $\cC_{X_0}(R)=\{gX_0g^{-1}, g\in
G(R)\}$ is given by a matrix in $\fg(R)$,
$$M:\begin{pmatrix}p_{ij}&q_{i\beta}\\r_{\alpha j
}&s_{\alpha\beta}\end{pmatrix}$$ whose entries satisfy
$p_1-c_1=0,\dots p_l-c_l=0$, and then  defines a homomorphism
$$M:A_{X_0}\longrightarrow R.$$ In view of the Observation
(\ref{obsfunctor}) we just have to show that for $f_1 \dots f_s \in R^o$,
$(f_1 \dots f_s)=(1)$, $F(R_{f_i})=\cC_{X_0}(R_{f_i})$.

Let $f\in R^\circ$. By the previous observation there  is an
obvious injective map $\cC_{X_0}(R_{f}) \rightarrow F(R_f)$. We
have to show that the $f_i$ can be chosen in such a way that the
map is also surjective.

This means that we need to prove that given $W \in \cC_{X_0}(R)$,
there exist $f \in R^\circ$ and $g \in G(R_{f})$ such that
\begin{equation}gWg^{-1}=D, \qquad\mathrm{or}\qquad gW=Dg\label{diagonal1}\end{equation}
 where $D$ is a diagonal matrix diagonal. Later
on we will prove that $D=X_0$.

We consider first a
superalgebra $R$ is a free superalgebra in the odd generators.
We decompose the matrices in  (\ref{diagonal1}) as sums of
matrices whose elements are homogeneous in the odd variables
\begin{equation*} (g_0+g_1+ \cdots )(W_0+W_1+ \cdots )= (D_0+D_1+ \cdots
)(g_0+g_1+ \cdots ).\end{equation*}
Then we can compare elements of the same degree obtaining
\begin{equation} g_0W_n+ g_1W_{n-1} + \cdots g_nW_0=D_0g_n+ D_1g_{n-1}+ \cdots
+D_ng_0. \label{diagonaln}\end{equation} We will prove the result
by induction. For $n=0$ we have \begin{equation} \label{diagonal0}
g_0W_0=D_0g_0
\end{equation}
But this is the classical result, with $D_0=X_0$. By the
hypothesis 2. of the theorem, $(p_1-c_1, \dots p_l-c_l)$ is the
ideal of the reduced orbit, so when we restrict to the category of
commutative algebras we know that $F$ is  the sheafification of
$\cC_{X_0}$ and so it  is represented by $\C[\fg^\circ]/(p_1-c_1
\dots p_l-c_l)$.

This means that there exist $f \in R^\circ$, $g_0 \in G(R^\circ)$,
such that $g_0W_0g_0^{-1}=D_0$. Moreover, the classical results
guarantees that  one can choose $f_1 \dots f_r$ among all possible
$f$'s in such a way that the ideal that they generate in
$R^\circ$, $(f_1 \dots f_r)=(1)=R^\circ$.

The induction proof is based on an argument given in Ref.
\cite{be}, page 117. We assume that the result is true up to order
$n-1$. Then we multiply (\ref{diagonaln}) by $g_0^{-1}$ to the
right. Using (\ref{diagonal0}) we obtain \begin{equation}
X_0g_ng_0^{-1}-g_ng_0^{-1}X_0+D_n=K_n,
\label{compute}\end{equation} where $K_n$ is a known matrix $$
K_n=-D_1g_{n-1}g_0^{-1}-D_2g_{n-2}g_0^{-1}- \dots +
g_0W_{n-1}g_0^{-1}+ \dots +g_{n-1}W_1g_0^{-1} $$
 and $D_n$ is a diagonal matrix. The matrix
$X_0g_ng_0^{-1}-g_ng_0^{-1}X_0$ has only elements outside the
diagonal, and if $(\lambda_i)$ are the eigenvalues of $X_0$ we
have that the entry $(ij)$ is given by
$$(X_0g_ng_0^{-1}-g_ng_0^{-1}X_0)_{ij}=(\lambda_i-\lambda_j)(g_ng_0^{-1})_{ij}.
$$ Then, $g_n$ and $D_n$  can be computed from (\ref{compute})
provided $\lambda_i\neq\lambda_j$.

To finish the proof we will have to show that $D_n=0$ for $n\geq
0$. Let us consider the case of $\fsl(m|n)$. Then the invariant
polynomials $p_i\in \C[\fsl(m|n)]$ that we have to consider are
$$p_i(M)=\str M^i, \qquad i=1,\dots m+n.$$ We want to prove that
for a diagonal matrix $D=D_0+D_1+\cdots$, if
\begin{equation}p_i(D)=p_i(D_0),\label{constraint}\end{equation} then $D=D_0$. Let
$$D=\mathrm{diag}(\lambda_1+{\lambda'}_1,\lambda_2+{\lambda'}_2,\dots
\lambda_{m+n}+{\lambda'}_{m+n}),$$ with ${\lambda'}_i$ contains
all the terms in the odd variables. Then (\ref{constraint})
implies the following homogeneous system: $$\begin{pmatrix}1&
\cdots &1&-1&\cdots&-1\\\lambda_1&\cdots
&\lambda_m&-\lambda_{m+1}&\cdots&-\lambda_{m+n}\\\cdots\\
\lambda_1^{m+n}&\cdots
&\lambda^{m+n}_m&-\lambda_{m+1}^{m+n}&\cdots&-\lambda^{m+n}_{m+n}
\end{pmatrix}\begin{pmatrix}\lambda'_1\\\lambda'_2\\\cdots\\\lambda'_{m+n}\\\end{pmatrix}$$
which can have non trivial solution solution
if the determinant $$\det\begin{pmatrix}1&
\cdots &1&-1&\cdots&-1\\\lambda_1&\cdots
&\lambda_m&-\lambda_{m+1}&\cdots&-\lambda_{m+n}\\\cdots\\
\lambda_1^{m+n}&\cdots
&\lambda^{m+n}_m&-\lambda_{m+1}^{m+n}&\cdots&-\lambda^{m+n}_{m+n}
\end{pmatrix}=(-1)^{nm}\prod_{i>j}(\lambda_i-\lambda_j)$$
is different from zero (Vandermonde determinant). So we have our
result if all the eigenvalues are different.

In the case of $\fosp_{m|n}$, the relevant polynomials are of the
form $\str M^{2i}$, so the result can be reproduced without
difficulty.

In the case we consider a superalgebra $R/J$ with $J$ an ideal, it
is enough to look at the images of the $(f_1,\cdots f_r)$ under
the projection $F/J$.

\hfill$\blacksquare$

\section{Deformation quantization of coadjoint orbits of
supergroups.\label{dco}}

Let $k=\C$. We start with the definitions of Poisson superalgebra
and its deformation.

\begin{definition}\label{poissonalgebra} Let $A$ be a commutative superalgebra. We say that $A$ is a
Poisson superalgebra if there exists a linear map (Poisson
superbracket)
\begin{eqnarray*}\{\;,\;\}:&A\otimes A\longrightarrow &A\\
&f\otimes g\longrightarrow &\{f,g\}\end{eqnarray*} such that

\smallskip

\noindent 1. $\{a,b\}=-(-1)^{p_ap_b}\{b,a\}$

 \noindent 2. $\{a,\{b,c\}\}+
 (-1)^{p_ap_b+p_ap_c}\{b,\{c,a\}\}+(-1)^{p_ap_c+p_bp_c}\{c,\{a,b\}\}=0$

 \noindent 3. $
\{a,bc\}=\{a,b\}c+(-1)^{p_ap_b}b\{a,c\} $

\smallskip
\noindent where $a,b,c$ are homogeneous elements of $A$ with
parities $p_a,p_b,p_c$\footnote{We want to remark that with this
definition we have only {\it even} Poisson
brackets.}.\end{definition} Let $\fg$ be a Lie superalgebra
with Lie superbracket $$[X_I,X_J]=c_{IJ}^KX_K$$ for a certain
homogeneous basis $\{X_I\}_{I=1}^{s+r}$ with the first $s$ vectors
even and the last $r$ odd. Then $\C[\fg^*]\simeq \C[x_1,\dots
x_s,x_{s+1},\cdots x_{s+r}]$ has a Poisson superalgebra structure
with superbracket given by
$$\{x_I,x_J\}=\sum_Kx_K([X_I,X_J])x_K=\sum_Kc_{IJ}^Kx_K$$ and
extended to the whole algebra by property 3. of Definition
\ref{poissonalgebra}.

Let $\fg$ be one of the Lie algebras considered in Theorem
\ref{orbittheorem}, which have a Cartan-Killing form, so
$\fg\simeq\fg^*$ and $\C[\fg]\simeq \C[\fg^*]$.  Let
$A=\C[\fg]/\cI$ be the algebra associated to the  superorbit of a
geometric element $X_0$ with the conditions of Theorem
\ref{orbittheorem}. Then $A$ is also a Poisson superalgebra with
the Poisson superbracket induced by the one in $\C[{\fg}]$. This
follows from the derivation property  3. in Definition
\ref{poissonalgebra} and the fact that $\hat p_i$ are invariant
polynomials.

In this section we want to construct a deformation quantization of
the superalgebra representing the orbit of a supergroup. We will
extend the method used in Refs. \cite{fl, fl2} for the classical
(non super) case.

\begin{definition} \label{deformation}Given a  Poisson superalgebra $A$,
a  formal
deformation (or deformation quantization) of $A$ is an associative
 non commutative superalgebra algebra $A_h$ over $\C[[h]]$, where $h$ is a formal
parameter, with the following properties:

\smallskip

\noindent 1. $A_h$ is isomorphic to $A[[h]]$ as a $\C[[h]]$-module.

\smallskip

\noindent 2. The multiplication $*_h$ in $A_h$ reduces mod($h$) to
the one in $A$.

\smallskip

\noindent 3. $\tilde a *_h \tilde b -\tilde  b *_h \tilde a =
h\{a,b\}$ mod $(h^2)$, where $\tilde a, \tilde b \in A_h$ reduce
to $a,b \in A$ mod($h$) and $\{\; ,\;\}$ is the Poisson
superbracket in $A$.
\end{definition}
Let $\fg$ be a Lie superalgebra and $\fg_h$ the Lie
superalgebra over $\C[[h]]$ obtained by multiplying the Lie bracket
of $\fg$ by the formal parameter $h$. Let us denote by $U_h$ the
universal enveloping algebra of $\fg_h$ (\cite{be} pg. 279).
 As in the classical case,
it is easy to prove that the associative, non commutative
superalgebra $U_h$ is a deformation quantization of $\C[\fg^*]$.

\bigskip

Let $\fg$ be a Lie superalgebra over $\C$ of the type considered
in Theorem \ref{orbittheorem}. By property 1 in Definition
\ref{deformation}, there exists an isomorphism of
$\C[[h]]$-modules $\psi:\C[\fg^*][[h]]]\rightarrow U_h$. We want
to prove that $\psi$ can be chosen in such a way that  there
exists an ideal $\cI_h\subset U_h$ such that
$\psi(\cI)=\cI_h$\footnote{the ideal $\cI$ is understood here as
an ideal of $\C[\fg^*][[h]]$}, and the map induced on the
quotients $\psi_h:\C[\fg^*]/\cI[[h]]\rightarrow U_{h}/I_{h}$ is an
isomorphism of  $\C[[h]]$-modules. We have then that the following
diagram

\begin{equation*}
\begin{CD}
\C[\fg^*][[h]]@>\psi>>U_{h}\\ @VV{\pi}V @VV{\pi_h}V\\
\C[\fg^*]/\cI[[h]]@>\psi_h>>U_{h}/\cI_{h}
\end{CD}
\label{cd}
\end{equation*}
with $\pi$ and $\pi_h$ the canonical projections, commutes.

\bigskip

We list first some known results about the enveloping algebra of
$\fg$, $U(\fg)$ \cite{be}. We denote by $S(\fg)$ the algebra of
super symmetric tensors on $\fg$. We have that
$S(\fg)\simeq\C[\fg^*]$. As before, $\{X_I\}_{I=1}^{s+r}$ denotes
a homogeneous basis of $\fg$ as vector
superspace. Let $\tau:\C[\fg^*]\rightarrow
U(\fg)$ be the supersymmetrizer map
\begin{equation*} \tau(x_{I_1}\cdots x_{I_p})=\frac{1}{p!}\sum_{s \in
S_p}(-1)^{\bar\sigma(s)}X_{s(I_1)}\otimes\cdots\otimes X_{s(I_p)}
\end{equation*}
where $S_p$ is the group of permutations of order $p$ and
$\bar\sigma(s)$ is the sign arising when performing the
permutation  $\begin{CD}X_{I_1}\otimes\cdots\otimes
X_{I_p}@>s>>X_{s(I_1)}\otimes\cdots\otimes X_{s(I_p)}\end{CD}$ as
if the homogeneous elements $X_{I_j}$ where supercommuting. $\tau$
is an isomorphism of $\fg$-modules.  Let $\C[\fg]^\fg$ denote the
set of polynomials invariant under the adjoint action of $\fg$ (in
particular $\hat p_i \in \C[\fg]^\fg$). Then  $\tau$ induces an
isomorphism of $\fg$-modules: $$\C[\fg]^{\fg} \simeq Z(U(\fg))$$
where $Z(U(\fg))$ is the center of $U(\fg)$

\bigskip

Let $\cI=(p_1-c_1, \dots p_l-c_l)$ as in Theorem
\ref{orbittheorem}. We set $$ \cI_h=(P_1-c_1(h) \dots P_l-c_l(h))
\subset U_h, \qquad P_i=\tau(\hat p_i), \quad c_k(0)=c_k, $$ with
$c_i(h)\in \C[[h]]$.

\begin{theorem} In the settings of Theorem
\ref{orbittheorem}, $U_h/\cI_h$ is a deformation quantization of
$\C[\fg^*]/\cI$. \end{theorem}

 {\bf Proof} The only property in
Definition \ref{deformation} which is not immediate is property
a., that is,  the fact of $U_h/\cI_h\simeq \C[\fg^*]/\cI[[h]]$ as
$\C[[h]]$-modules. The proof can be almost translated from  the
classical case in Ref. \cite{fl}. Let $\{x_{I_1}, \dots ,
x_{I_k}\}_{(I_1,\dots , I_k) \in S}$ be a basis of $\C[\fg^*]/I$
as $\C$-module, where $S$ is a fixed set of multiindices
appropriate to describe the basis. In particular, we can take them
such that $I_1 \leq \cdots\leq I_k$. Proving that the monomials
$\mathcal{B}=\{X_{I_1}\cdots X_{I_k}\}_{(I_1, \dots , I_k) \in S}$
are a basis for  $U_h/I_h$ will be enough. The proof that
$\mathcal{B}$ is a generating set for $U_h/I_h$ is identical to
the proof of Proposition 3.13 in Ref. \cite{fl} and we will not
repeat it here. For the linear independence, we have to show that
there is no  relation among them. Suppose that $G\in \cI_h$ is
such relation, $$G=G_0+G_1h+ \cdots , \qquad G_i \in
\hbox{span}_{\C}\{X_{I_1}\cdots X_{I_k}\}_{(I_1 \dots I_k)\in S}.
$$ Assume $G_i=0$,  $i<r$, $G_r \neq 0$ so we can write $G=h^rF$,
with
\begin{equation} F=F_0+hF_1+ \cdots , \qquad F_0 \neq 0.
\label{tobereduced}\end{equation}

We need to prove the non trivial fact that if  $hF \in \cI_h$
 then $F \in \cI_h$.  We will denote by a capital letter, say $P$, an element in $U_h$,
 by $\hat p$ its  projection onto $\C[\fg^*]$ and by $p$ its further projection onto $\C[{\fg^\circ}^*]$.
 Assume that
 \begin{equation}hF=\sum_iA_i(P_i-c_i(h))\label{torsion}\end{equation}
with $c_i(h)=c_i+c_i^1h+ \dots c_nh^n$.

Then, reducing mod $h$ (taking $h=0$) we obtain
 \begin{equation}0=\sum_i\hat a_i(\hat p_i-c_i).\label{reduced}\end{equation}
 Setting all the odd variables to 0, we have that (\ref{reduced}) implies that there exist
 $b_{ij}$, antisymmetric in $i$ and $j$ such that
 $$a_i=\sum_jb_{ij}(p_j-c_j)$$
(see for example Ref. \cite{ol}, pg 81, or lemma 3.8 of Ref.
\cite{fl}), provided that the differentials $dp_i$ are independent
on the orbit (condition which is satisfied \cite{va} in our case).
The generalization of this property to the supersymmetric case
deserves special treatment, but the proof is rather technical and
we will do it separately in Lemma \ref{lemmaolver}. Assuming that
this is true, then there exist $\hat b_{ij}$, antisymmetric in $i$
and $j$ such that $$\hat a_i=\sum_j\hat b_{ij}(\hat p_j-c_j).$$ It
is easy to convince oneself that this equation can be lifted to
$U_h$, so there exists $A'_i$ and $B_{ij}$ antisymmetric in $i$
and $j$ such that $$A'_i=\sum_j B_{ij}(P_j-c_j(h)),$$ with
$A_i=A_i'$ mod $h$, since they both project to $a_i$, i.e.
$A_i=A_i'+hC_i$. If one substitutes in (\ref{torsion}):
$$
hF=h\sum_iC_i(P_i-c_i(h))
$$
since $\sum_j B_{ij}(P_j-c_j(h))(P_i-c_i(h))=0$ (the $B_{ij}$'s
are antisymmetric). Hence we get the fact: $hF \in \cI_h$
then $F \in \cI_h$.

So being $F\in\cI_h$, we can reduce (\ref{tobereduced}) mod $h$, $$ f = \sum
a_i(\hat p_i-c_{i}). $$ But $f$ would represent a non trivial
relation among the monomials \break $\{x_{I_1}\cdots
x_{I_k}\}_{(I_1 ... I_k) \in S}$ in $\C[\fg^*]/\cI$, which is a
contradiction, so the linear independence is proven.
\hfill$\blacksquare$

\begin{lemma}\label{lemmaolver}. Let $A$ be the free commutative superalgebra over $\C$ generated by
$M$ even variables $x_1 \dots x_M$ and $N$ odd variables $\xi_1
\dots \xi_N$. Let $ q_1 \dots  q_n$ even polynomials in $A$ and
denote by $q^\circ_1,\dots q^\circ_n$ their projections onto
$A^\circ$. Assume that the  $q^\circ_i$'s satisfy the following
property:

\medskip

\noindent If $\sum_i q^\circ_if^\circ_i=0$ for some $f^\circ_i\in
A^\circ$, then there exist $F^\circ_{ij} \in A^\circ$ such that
$f^\circ_i=\sum_{j} F^\circ_{ij}q^\circ_j$, with
$F^\circ_{ij}=-F^\circ_{ji}$, $i,j=1 \dots n$.

\medskip

\noindent Then, if $\sum_i q_i f_i=0$ for some $f_i\in A$, there
exist $ F_{ij} \in A$ such that $ f_i=\sum_{j}  F_{ij} q_j$, with
$ F_{ij}=- F_{ji}$, $i,j=1 \dots n$.
\end{lemma}

{\bf Proof}.
We write
\begin{eqnarray*}&&f_i=\sum_{k=0}^N\sum_{1\leq \alpha_1<\dots
    \alpha_k\leq N}f_i^{\alpha_1\dots
    \alpha_k}\xi_{\alpha_1}\cdots\xi_{\alpha_k}\\&&q_i=
\sum_{l=0}^N\sum_{1\leq \beta_1<\dots
    \beta_l\leq N}q_i^{\beta_1\dots
    \beta_l}\xi_{\beta_1}\cdots\xi_{\beta_l}.\end{eqnarray*} Since $q_i$ are even functions we have that
$l$ is always even number (that is, the components with $l$ odd are zero). The terms for $k=0$ or $l=0$ correspond to
$f_i^\circ$ and $q_i^\circ$ respectively.

The condition $\sum_{i=1}^nf_iq_i=0$  reads in components
\begin{equation}\sum_{i=1}^n\sum_{s=0}^p \sum_{1\leq k_1<\cdots < k_s\leq p}(-1)^{\cl(k_1,\dots k_s,1,\dots \hat k_1,\dots \hat k_s,\dots p)}
 f_i^{\alpha_{k_1}\cdots \alpha_{k_s}}q_i
^{\alpha_1\cdots \hat \alpha_{k_1}\cdots \hat\alpha_{k_s}\cdots \alpha_p}=0,\label{hypothesis}\end{equation}
for $p=1, \dots N$ and where $\cl(k_1,\dots k_s,1,\dots \hat k_1,\dots \hat k_s,\dots p)$ is the length of the permutation
 $(k_1,\dots k_s,1,\dots \hat k_1,\dots \hat k_s,\dots p)$.  For $p=0$ we have
\begin{equation}\sum_{i=1}^n f_i^\circ q_i^\circ=0.\label{hypo0}\end{equation}

\bigskip

We will show the lemma by induction. Let $\cI^\mathrm{odd}_m$ be the
ideal in $A$ generated by the products
$$\cI^\mathrm{odd}_m=(\xi_{\alpha_1}\cdots \xi_{\alpha_{m+1}}, \dots,
\xi_1\cdots \xi_N),$$
and consider the projections
$\pi_m:A\rightarrow A/\cI^\mathrm{odd}_m$. Assume that
there exists
$$F_{ij}^m=\sum_{r=0}^mF_{ij}^{\delta_1\cdots
  \delta_r}\xi_{\delta_1}\cdots \xi_{\delta_r},$$
with $F_{ij}^m=-F^m_{ji}$, such that
\begin{equation}\pi_s(f_i)=\pi_s(F_{ij}^mp_j)\label{inductionhypothesis}\end{equation}
for all $s\leq m$. (To simplify the notation the sum over
reapeated indices $i$ and $j$ will be understood).  For $m=0$ this is the classical condition, that is
guaranteed by (\ref{hypo0}) and the  hypothesis of the theorem. We must show that there exist
$F_{ij}^{\delta_1\cdots \delta_{m+1}}$ such that

$$F_{ij}^{m+1}=\sum_{r=0}^{m+1}F_{ij}^{\delta_1\cdots
  \delta_r}\xi_{\delta_1}\cdots \xi_{\delta_r}$$  satisfies
$$\pi_s(f_i)=\pi_s(F_{ij}^{m+1}p_j)$$
for all $s\leq m+1$. Since $\pi_N$ is the identity, we will have our
result.

\bigskip

Let us write the induction hypothesis (\ref{inductionhypothesis}) in
components,
\begin{equation}f_i^{\alpha_1\cdots \alpha_s}=\sum_{r=0}^s\sum_{1\leq l_1<\cdots l_r\leq s}
(-1)^{\cl(l_1\dots l_r,1,\dots \hat l_1,\dots \hat l_r,\dots s)}
F_{ij}^{\alpha_{l_1}\dots \alpha_{l_r}}q_j^{\alpha_1\dots\hat\alpha_{l_1}\dots \hat\alpha_{l_r},\dots \alpha_s},
\label{inductionhypothesis2}\end{equation}
with $s\leq m$.
We substitute this expression for $s\leq m$ in (\ref{hypothesis}). Taking
$p=m+1$  we have

\begin{eqnarray}
&&f_i^{\alpha_1\dots
 \alpha_{m+1}}q_i^0+\sum_{s=0}^m\sum_{1\leq k_1<\cdots <k_s\leq m+1}
\sum_{r=0}^s\sum_{k_1\leq k_{l_1}<\cdots <k_{l_r}\leq k_s}\nonumber\\
&&(-1)^{\cl(k_1,\dots
    k_s,1,\dots \hat k_1,\dots \hat k_s,\dots m+1)}
(-1)^{\cl(k_{l_1}\dots k_{l_r}k_1\dots \hat k_{l_1}\dots \hat k_{l_r}
\dots k_s)}\nonumber\\
&&F_{ij}^{\alpha_{k_{l_1}}\dots \alpha_{k_{l_r}}}q_j^{\alpha_{k_1}\dots\hat\alpha_{k_{l_1}}\dots \hat\alpha_{k_{l_r}}
 \dots \alpha_{k_{s}}}
q_i^{\alpha_1,\dots \hat\alpha_{k_1}\dots \hat \alpha_{k_s}\dots
  \alpha_{m+1}}=0\label{allterms}\end{eqnarray}

We distinguish two kinds of terms in (\ref{allterms}), the ones that are multiplied by
$q_i^\circ$ and the ones that are not. A generic term that does not contain $q_i^\circ$ is of the form

\begin{eqnarray}
(-1)^{\cl(k_1,\dots k_s,1,\dots \hat k_1,\dots \hat k_s,\dots m+1)}
(-1)^{\cl(k_{l_1}\dots k_{l_r}k_1\dots \hat k_{l_1}\dots \hat k_{l_r}\dots k_s)}\nonumber\\
F_{ij}^{\alpha_{k_{l_1}}\dots \alpha_{k_{l_r}}}q_j^{\alpha_{k_1}\dots\hat\alpha_{k_{l_1}}\dots \hat\alpha_{k_{l_r}}
 \dots \alpha_{k_{s}}}
q_i^{\alpha_1,\dots \hat\alpha_{k_1}\dots \hat \alpha_{k_s}\dots \alpha_{m+1}}\label{term}
\end{eqnarray}
with $r\neq s$ and $s\neq m+1$.
We are going to see that these terms cancel. The reason is that  for each term as (\ref{term}) there exists an identical
 term where $q_i$ is interchanged with $q_j$. Then, they cancel because of the antisymmetry of $F_{ij}$.
 We first consider  the following   set of indices
$$\{k'_1 , \dots k'_{s'}\}=\bigr(\{1, \dots m+1\}-\{k_1 , \dots
 k_{s}\}\bigl)\cup \{k_{l_1}, \dots k_{l_r}\}.$$
Then
\begin{eqnarray}\{k'_1 , \dots k_{s'}'\}- \{k_{l_1}, \dots k_{l_r}\}=\{1, \dots m+1\}-\{k_1 , \dots
 k_{s}\}\nonumber\\
\{k_1 , \dots k_{s}\}- \{k_{l_1}, \dots k_{l_r}\}=\{1, \dots m+1\}-\{k'_1 , \dots
 k_{s'}'\}\label{permutations}
\end{eqnarray}
 and $s'=m+1-s+r$. By construction $\{k_{l_1}, \dots k_{l_r}\}\subset
 \{k'_1 , \dots k'_{s'}\}$, so we can consider the term
 \begin{eqnarray}
(-1)^{\cl(k'_1,\dots k'_{s'},1,\dots \hat k'_1,\dots \hat k'_{s'},\dots m+1)}
(-1)^{\cl(k_{l_1}\dots k_{l_r}k'_1\dots \hat k_{l_1}\dots \hat k_{l_r}\dots k'_{s'})}\nonumber\\
F_{ij}^{\alpha_{k_{l_1}}\dots \alpha_{k_{l_r}}}q_j^{\alpha_{k'_1}\dots\hat\alpha_{k_{l_1}}\dots \hat\alpha_{k_{l_r}}
 \dots \alpha_{k_{s'}'}}
q_i^{\alpha_1,\dots \hat\alpha_{k'_1}\dots \hat \alpha_{k_{s'}'}\dots \alpha_{m+1}}\label{term2}
\end{eqnarray}
(\ref{term2}) is equal to (\ref{term}) with $q_i$ in the place of
$q_j$, except possibly by a sign. We are going to see that the sign is
the same. We start by observing that

\begin{eqnarray*}
&&\cl(k_1 \dots k_s, 1 \dots \hat k_1 \dots \hat k_s \dots m+1)=
\cl(k_{l_1} \dots k_{l_r},1 \dots \hat k_1 \dots \hat k_s \dots m+1)+\\
&&\cl(k_1 \dots \hat k_{ l_1} \dots \hat k_{l_r} \dots k_s,
1 \dots \hat k_1 \dots \hat k_s \dots m+1),
\\
\\
&&\cl(k_1' \dots k'_{s'}, 1 \dots \hat k_1' \dots \hat k_{s'}' \dots
m+1)=\cl(k_{l_1} \dots k_{l_r},1 \dots \hat k_1' \dots \hat k'_{s'} \dots m+1)+
\\
&&\cl(k_1' \dots \hat k_{l_1} \dots \hat k_{ l_r} \dots k'_{s'},
1 \dots \hat k_1' \dots \hat k'_{s'} \dots m+1).
\end{eqnarray*}
It is now crucial the fact that $s-r$, $m+1-s$, $s'-r$, $m+1-s'$ are
all even numbers.  The last two terms in the above equalities are
easily seen equal modulo 2 by using (\ref{permutations}). Then we have

\begin{eqnarray*}
&&\cl(k_1 \dots k_s, 1 \dots \hat k_1 \dots \hat k_s \dots m+1)-
\cl(k_1' \dots k'_{s'}, 1 \dots \hat k_1' \dots \hat k_{s'}' \dots
m+1)=\\
&&\cl(k_{l_1} \dots k_{l_r},1 \dots \hat k_1 \dots \hat k_s \dots m+1)
-\cl(k_{l_1}
 \dots k_{l_r},1 \dots \hat k_1' \dots \hat k'_{s'} \dots
 m+1)\;\mathrm{mod}2=\\
&&\cl(k_{l_1} \dots k_{l_r}, k_1', \dots \hat k_{l_1} \dots \hat k_{l_r},\dots k_{s'}' )
-\cl(k_{l_1}
 \dots k_{l_r},  k_1, \dots \hat k_{l_1} \dots \hat k_{l_r},\dots  k_{s})\;\mathrm{mod}2,
\end{eqnarray*}
where we have used (\ref{permutations}) again. We have proven our
claim and all the terms in (\ref{allterms}) that do not contain
$q_i^0$ cancel. Then from (\ref{allterms}) we deduce

\begin{eqnarray*}
&&\bigl(f_i^{\alpha_1\dots
 \alpha_{m+1}}-\sum_{s=0}^m\sum_{1\leq k_1<\cdots <k_s\leq m+1}
(-1)^{\cl(k_1,\dots
    k_s,1,\dots \hat k_1,\dots \hat k_s,\dots m+1)}
\nonumber\\
&&F_{ij}^{\alpha_{k_{1}}\dots \alpha_{k_s}}
q_j^{\alpha_1,\dots \hat\alpha_{k_1}\dots \hat \alpha_{k_s}\dots
  \alpha_{m+1}}\bigr)  q_i^{\circ}=0.\end{eqnarray*}

From the hypothesis of the theorem, this means that there exists
 $F_{ij}^{\alpha_1\dots \alpha_{m+1}}$, antisymmetric in $i$ and $j$
 such that
\begin{eqnarray*}
&&f_i^{\alpha_1\dots
 \alpha_{m+1}}-\sum_{s=0}^m\sum_{1\leq k_1<\cdots <k_s\leq m+1}
(-1)^{\cl(k_1,\dots
    k_s,1,\dots \hat k_1,\dots \hat k_s,\dots m+1)}
\nonumber\\
&&F_{ij}^{\alpha_{k_{1}}\dots \alpha_{k_s}}
q_j^{\alpha_1,\dots \hat\alpha_{k_1}\dots \hat \alpha_{k_s}\dots
  \alpha_{m+1}}=F_{ij}^{\alpha_1\dots
 \alpha_{m+1}}q_j^\circ,\end{eqnarray*}
from which one can compute
$$f_i^{\alpha_1\dots
 \alpha_{m+1}}=\sum_{s=0}^{m+1}\sum_{1\leq k_1<\cdots <k_s\leq m+1}
(-1)^{\cl(k_1,\dots
    k_s,1,\dots \hat k_1,\dots \hat k_s,\dots m+1)}F_{ij}^{\alpha_{k_{1}}\dots \alpha_{k_s}}
q_j^{\alpha_1,\dots \hat\alpha_{k_1}\dots \hat \alpha_{k_s}\dots
  \alpha_{m+1}},$$
terminating our induction.

  \hfill$\blacksquare$

  \section*{Acknowledgments}

  We would like to thank Prof.  V. S. Varadarajan and Prof. A. Vistoli for many
helpful comments.

 Work supported in part by  the European
Community's Human Potential Program under contract
HPRN-CT-2000-00131 Quantum Space-Time, M. A. Ll. is associated to
Torino University.

The work of M. A Ll. has also been supported by the research grant
BFM 2002-03681 from the Ministerio de Ciencia y Tecnolog\'{\i}a
(Spain) and from EU FEDER funds.

\end{document}